\documentclass[reqno,10pt]{amsart}
\usepackage{mathrsfs,graphicx,amsthm,amsfonts,amsxtra,xspace,array,cite,epic,float,ifpdf}
\usepackage[left=2.5cm]{geometry}
\newcommand\version{public}
\newcommand\finalized{yes}
\usepackage{fullpage,setspace}
\usepackage{tfrupee}  

\usepackage{mathrsfs} 
\newcommand\choosefont[1]{\usepackage{#1}}
\usepackage[bottom]{footmisc}  
\usepackage{enumerate}   
\usepackage{url,cite,fancyheadings}
\usepackage{enumerate}
\usepackage{asymptote}
\usepackage{amsmath,amssymb,amscd,mathrsfs,latexsym}
		\usepackage{rsfso}  
		\usepackage{mathtools}
\usepackage{mathdots}

\usepackage{tikz}

\usepackage{ifthen}     

\newcommand\pubpri[2]{%
\ifthenelse{\equal{\version}{public}}%
{{#1}}%
{\ifthenelse{\equal{\finalized}{no}}{\marginpar{\scshape\small Pubpri Alert}{#2}}{#2}}{}}
\newcommand\pubprinoalert[2]{%
\ifthenelse{\equal{\version}{public}}%
{{#1}}%
{#2}}

\newcommand\ignore[1]{}

\usepackage{color}
\providecommand\wantcolor{yes}   %
\definecolor{backgroundyellow}{cmyk}{.2,.1,.8,.2}
\definecolor{backgroundblue}{rgb}{0,0,1}
\definecolor{backgroundred}{rgb}{1,0,0}
\definecolor{backgroundmagenta}{cmyk}{0,1,0,0}
\ifthenelse{\equal{\wantcolor}{no}}{\renewcommand\color[1]{}}{}

\usepackage{enumitem}

\newcommand\mysection{\section}

\newcommand\mysubsection{\subsection}

\newcommand\mysubsubsection[1]{%
		\subsubsection{\sffamily\upshape\mdseries #1}}

\newcommand\mysss{\mysubsubsection}

\usepackage{chngcntr}
\counterwithin{figure}{section}

\usepackage{url,hyperref}  
\hypersetup{colorlinks=true,
linkcolor=blue,
filecolor=magenta,
urlcolor=cyan,}
\usepackage{graphicx,wrapfig,epstopdf}

\providecommand{\theoremnumbering}{document}

\ifthenelse{\equal{\theoremnumbering}{section}}{\newtheorem{annotation}{Annotation}[section]}%
{\ifthenelse{\equal{\theoremnumbering}{subsection}}{\newtheorem{annotation}{Annotation}[subsection]}{\newtheorem{annotation}{Annotation}}}
\newtheorem{theorem}[annotation]{
		Theorem}
\newtheorem{lemma}[annotation]{
		Lemma}
\newtheorem{definition}[annotation]{
		Definition}
\newtheorem{corollary}[annotation]{
		Corollary}
\newtheorem{proposition}[annotation]{
		Proposition}

\newtheorem{example}[annotation]{
		Example}
\newcommand\bexample{\begin{example}\begin{rm}}
\newcommand\eexample{\end{rm}\hfill$\Box$\end{example}}
\newtheorem{examplenobox}[annotation]{
		Example}
\newcommand\bexamplenobox{\begin{examplenobox}\begin{rm}}
\newcommand\eexamplenobox{\end{rm}\end{examplenobox}}
\newtheorem{exercise}[annotation]{
		Exercise}
\newcommand\bexercise{\noindent\begin{exercise}\begin{rm}}
\newcommand\eexercise{\end{rm}\end{exercise}}
\newtheorem{notation}[annotation]{
		Notation}
\newcommand\bnotation{\begin{notation}\begin{rm}}
\newcommand\enotation{\end{rm}\end{notation}}

\newtheorem{remark}[annotation]{
		Remark}
\newcommand\bremark{\begin{remark}
\begin{upshape}}
\newcommand\eremark{\end{upshape}
\end{remark}}
\newenvironment{remark*}{%
\par\noindent{\scshape 
  Remark: }\begin{rm}}{\hfill\end{rm}\newline} 
\newcommand\bremarkstar{\begin{remark*}}
\newcommand\eremarkstar{\end{remark*}}
\newcommand\bdefn{\begin{definition}
\begin{upshape}}
\newcommand\edefn{\end{upshape}
\end{definition}}
\newtheorem{caveat}[annotation]{
		Caveat}
\newcommand\bcaveat{\begin{caveat}
\begin{upshape}}
\newcommand\ecaveat{\end{upshape}
\end{caveat}}
\newenvironment{caveatstar}{
\par\noindent{\scshape\bfseries
  Caveat: }\begin{rm}}{\end{rm}\newline} 
\newcommand\bcaveatstar{\begin{caveatstar}}
\newcommand\ecaveatstar{\end{caveatstar}}
\newenvironment{myproof}{%
\par\noindent{\scshape 
  Proof: }\begin{rm}}{\hfill$\Box$\end{rm}\newline} 
\newcommand\bmyproof{\begin{myproof}}
\newcommand\emyproof{\end{myproof}}
\newenvironment{myproofnobox}{%
\par\noindent{\scshape Proof: }\begin{rm}}{\end{rm}\hfill\newline}
\newcommand\bmyproofnobox{\begin{myproofnobox}}
\newcommand\emyproofnobox{\end{myproofnobox}}
\newenvironment{myproofof}[1]{%
\par\noindent{\scshape 
  Proof (of~#1): }\begin{rm}}{\hfill$\Box$\end{rm}\newline} 
\newcommand\bmyproofof{\begin{myproofof}}
\newcommand\emyproofof{\end{myproofof}}
\newenvironment{myproofofnobox}[1]{%
\par\noindent{\scshape 
  Proof (of~#1): }\begin{rm}}{\end{rm}\hfill\newline} 
\newcommand\bmyproofofnobox{\begin{myproofofnobox}}
\newcommand\emyproofofnobox{\end{myproofofnobox}}
\newenvironment{solution}{%
\par\noindent{\scshape Solution: }\begin{rm}}{\hfill$\Box$\end{rm}\newline}
\newenvironment{solutionnobox}{%
\par\noindent{\scshape Solution: }\begin{rm}}{\end{rm}}
\newcommand\bsolution{\begin{solution}\begin{rm}}
\newcommand\esolution{\end{rm}\end{solution}}
\newcommand\bsolutionnobox{\begin{solutionnobox}\begin{rm}}
\newcommand\esolutionnobox{\end{rm}\end{solutionnobox}}
\newcommand\bthm{\begin{theorem}}
\newcommand\ethm{\end{theorem}}
\newcommand\bcor{\begin{corollary}}
\newcommand\ecor{\end{corollary}}
\newcommand\blemma{\begin{lemma}}
\newcommand\elemma{\end{lemma}}
\newcommand\bprop{\begin{proposition}}
\newcommand\eprop{\end{proposition}}
\newcommand\beqn{\begin{equation}}
\newcommand\eeqn{\end{equation}}
\newcommand\beqnstar{\begin{equation*}}
\newcommand\eeqnstar{\end{equation*}}
\usepackage{amsthm}
\usepackage{textcomp}
\usepackage{calrsfs,mathrsfs}  

\newcommand\mtitle[1]%
{\marginpar{\sffamily\raggedright\upshape\small #1}}

\providecommand\finalized{no}
\ifthenelse{\equal{\finalized}{yes}}{}{\marginparwidth=2cm}
\ifthenelse{\equal{\finalized}{yes}}%
{\newcommand\mylabel[1]{\label{#1}}}%
{\newcommand\mylabel[1]{\label{#1}\marginpar{[{\ttfamily\upshape\tiny #1}]}}}
\ifthenelse{\equal{\finalized}{yes}}%
{\newcommand\checked[1]{}}%
{\newcommand\checked[1]{\marginpar{[{\ttfamily\upshape\tiny CHECKED: #1}]}}}
\ifthenelse{\equal{\finalized}{yes}}%
{\newcommand\spellchecked[1]{}}%
{\newcommand\spellchecked[1]{\marginpar{[{\ttfamily\upshape\tiny SPELLCHECKED: #1}]}}}

\providecommand\version{public}   
\ifthenelse{\equal{\version}{public}}%
{\newcommand\mcomment[1]{}}%
{\newcommand\mcomment[1]{\marginpar{{\raggedright\sffamily\upshape\small
\begin{spacing}{0.75} #1\end{spacing}}}}}
\ifthenelse{\equal{\version}{public}}%
{\newcommand\fcomment[1]{}}%
{\newcommand\fcomment[1]{\footnote{#1}}}
\ifthenelse{\equal{\version}{public}}%
{\newcommand\comment[1]{}}%
{\newcommand\comment[1]{{\small #1}}}

\newcommand\st{\,|\,}

\newcommand\Vlambda{V(\lambda)}
\newcommand\vlambda{v_\lambda}
\newcommand\pseq{\underline{\pi}}
\newcommand\patl{\pattern(\lseq)}
\newcommand\patternp{\pattern'}
\newcommand\gtpatternp{\gt_{\patternp}}
\newcommand\patternt\tpattern
\newcommand\patternq{\mathcal{Q}}
\newcommand\gtpatternt{\gt_{\patternt}}
\newcommand\cprime{$'$}
\newcommand\germ{\mathfrak}

\newcommand{\be}{\begin{enumerate}}
\newcommand{\ee}{\end{enumerate}}
\newcommand{\beq}{\begin{equation}}
\newcommand{\eeq}{\end{equation}}

\renewcommand\omega{\varpi}

\newcommand\lieg{\mathfrak{g}}

\providecommand\lieg{\mathfrak{g}}
\providecommand\lieh{\mathfrak{h}}
\providecommand\lieb{\mathfrak{b}}

\providecommand\lseq{\underline{\lambda}}

\providecommand\lpseq{\underline{\lambda}'}

\providecommand\mseq{\underline{\mu}}

\providecommand\vcl[1]{v_{#1}}
\providecommand\rhocl[1]{\Theta_{#1}}

\DeclareMathOperator{\weight}{wt}

\providecommand\slrpone{\mathfrak{sl}_{r+1}}
\providecommand\currlieg{\lieg[t]}
\providecommand\curralg\currlieg

\providecommand\germ\mathfrak

\providecommand\lseqrpone{\lseq^{r+1}}
\providecommand\lseqone{\lseq^1}
\providecommand\lseqr{\lseq^r}

\providecommand\lpseqone{\lpseq^1}
\providecommand\lpseqr{\lpseq^r}
\providecommand\lpseqrpone{\lpseq^{r+1}}

\providecommand\pattern{\mathcal{P}}

\providecommand\beq{\begin{equation}}
\providecommand\eeq{\end{equation}}

\providecommand\cprime{$'$}

\providecommand{\be}{\begin{enumerate}}
\providecommand{\ee}{\end{enumerate}}

\renewcommand\omega{\varpi}

\providecommand\gt{\xi}
\providecommand\gtpattern{\gt_\pattern}
\providecommand\ellki{\ell^k_i}

\providecommand\coeff{c}
\DeclareMathOperator{\lengthp}{\textup{length}}
\providecommand\tpattern{\tilde{\pattern}}
\providecommand\tell{\tilde{\ell}}
\providecommand\tlambda{\tilde{\lambda}}
\providecommand\tlseq{\tilde{\lseq}}
\providecommand\ltildeseq{\tilde{\lseq}}

\numberwithin{equation}{section} 
\title[GT bases versus CL bases]{A relationship between Gelfand-Tsetlin bases and Chari-Loktev bases \\
for irreducible finite dimensional representations\\ of special linear Lie algebras}
\author{K.~N.~Raghavan}
\address{The Institute of Mathematical Sciences, HBNI\\ CIT campus, Taramani\\ Chennai 600113, India}
\email{knr@imsc.res.in}
\author{B.~Ravinder}
\address{Chennai Mathematical Institute,  H1 Sipcot IT Park, Kelambakkam, Siruseri, Chennai 603103, INDIA}
\email{bravinder@cmi.ac.in}
\author{Sankaran Viswanath}
\address{The Institute of Mathematical Sciences, HBNI\\ CIT campus, Taramani\\ Chennai 600113, India}
\email{svis@imsc.res.in}
\thanks{\noindent KNR and SV acknowledge support under a DAE project grant.
        BR acknowledges support from DST under the INSPIRE Faculty scheme (DST/INSPIRE/04/2016/001471)
and from the Infosys Foundation.}

\choosefont{newcent}
	\ignore{
	\begingroup\makeatletter\ifx\SetFigFont\undefined%
	\gdef\SetFigFont#1#2#3#4#5{%
  	\reset@font\fontsize{#1}{#2pt}%
  	\fontfamily{#3}\fontseries{#4}\fontshape{#5}%
  	\selectfont}%
	\fi\endgroup%
	}
\begin{document}
\allowdisplaybreaks
\begin{abstract}
	We consider two bases for an arbitrary finite dimensional irreducible representation of a complex special linear Lie algebra: the classical Gelfand-Tsetlin basis and the relatively new Chari-Loktev basis. Both are parametrized by the set of (integral Gelfand-Tsetlin) patterns with a fixed bounding sequence determined by the highest weight of the representation.  We define the {\em row-wise dominance} partial order on this set of patterns, and  prove that the transition matrix between the two bases is triangular with respect to this partial order.  We write down explicit expressions for the diagonal elements of the transition matrix.
\end{abstract}

\maketitle
\mysection{Introduction}\mylabel{s:introduction}
\noindent
Let $\lieg=\slrpone$ be the simple Lie algebra of complex matrices of size $(r+1)\times(r+1)$ and trace zero.
There are various bases known for (finite dimensional) irreducible representations of~$\lieg$.  
Here we consider two of them.   The first is the well known basis constructed by Gelfand-Tsetlin in~\cite{GT}.  We consider it as reformulated by Molev in~\cite[Theorem~4.3]{molevhb}.  The details are recalled in~\S\ref{s:gt} below.

The second basis is of more recent vintage and is obtained as a by-product of the construction by Chari-Loktev~\cite{cladv2006} of bases for local Weyl modules of the current algebra of~$\lieg$. 
We consider the construction of Chari-Loktev as reformulated by the present authors in~\cite[Theorem~4.3]{rrv:areamax}.  Dominant integral weights of~$\lieg$ form a natural parametrizing set for local Weyl modules of the current algebra (as they do for irreducible representations of~$\lieg$).    Local Weyl modules admit a natural grading by non-negative integers, and the grade zero slice of a local Weyl module is precisely the corresponding irreducible representation of~$\lieg$.  The Chari-Loktev bases for local Weyl modules are graded.   Thus,  by choosing only those elements that have grade zero from the Chari-Loktev basis for a local Weyl module,   we obtain a basis for the corresponding irreducible representation of~$\lieg$.   The details are recalled in~\S\ref{s:cl} below.

Both the Gelfand-Tsetlin and the Chari-Loktev bases (for an irreducible representation of~$\lieg$) are parametrized by (integral Gelfand-Tsetlin) patterns with a fixed bounding sequence determined by the highest weight of the representation.   The requisite details about patterns are recalled in~\S\ref{s:patterns} below.   It is natural to wonder if there is a relationship between the two bases.  The goal of the present note is Theorem~\ref{t:triangular} in~\S\ref{s:triangular} below,  which shows that the transition matrix between the two bases is triangular with respect to the {\em row-dominance\/} partial order (defined in the beginning of~\S\ref{s:triangular}) on their common parametrizing set.   Equation~\eqref{e:tricoeff} gives an explicit formula for the diagonal entries of the transition matrix.  A simple example illustrating Theorem~\ref{t:triangular} is given in~\S\ref{s:example}.

\mysection{Notation and terminology}\mylabel{s:notation}
\noindent 
The following notation and terminology remain fixed throughout this note.    Let~$\lieg=\slrpone$ be the simple Lie algebra of complex matrices of size $(r+1)\times(r+1)$ and trace zero.   Let $\lieh$ be the Cartan subalgebra of~$\lieg$ consisting of its diagonal matrices,  and $\lieb$ the Borel subalgebra of $\lieg$ consisting of its upper triangular matrices.    When we speak of roots, weights, dominant weights, etc.\ in the sequel,   it is with respect to this choice of the Cartan and Borel subalgebras.

Throughout, $\lambda$ denotes a fixed dominant integral weight,
$V(\lambda)$ the irreducible
representation of $\lieg$ with highest weight~$\lambda$, and
$v_\lambda$ a non-zero element in the $\lambda$-weight space of
$V(\lambda)$.  Note that $v_\lambda$ is unique up to scaling.  

We denote by $\lseq$ an arbitrarily fixed member of the equivalence class of $(r+1)$-tuples associated to $\lambda$ as in~\S\ref{ss:tuples} below.  The results below are so formulated that the choice of~$\lseq$ does not matter.
\mysubsection{Weights and associated tuples}\mylabel{ss:tuples}
\noindent
Recall that a {\em weight\/} is any linear functional on~$\lieh$.    For $i$, $1\leq i\leq r+1$, let $\epsilon_i$ denote the weight that maps an element of $\lieh$ to its diagonal entry in position $(i,i)$.  Evidently, the $\epsilon_i$, $1\leq i\leq r+1$, span the space of all weights, so any weight can be written as~$a_1\epsilon_1+\cdots+a_{r+1}\epsilon_{r+1}$ for some complex numbers $a_1$, \ldots, $a_{r+1}$.  Given an $(r+1)$-tuple $(a_1,\ldots,a_{r+1})$ of complex numbers,   we associate with it the weight $a_1\epsilon_1+\cdots+a_{r+1}\epsilon_{r+1}$.  Since the trace of any element of $\lieh$ is zero, two such tuples differ from each other by a scalar multiple of the constant $(r+1)$-tuple $(1,\cdots,1)$ if and only if the weights associated to them are the same.    

A weight $a_1\epsilon_1+\cdots+a_{r+1}\epsilon_{r+1}$ is {\em integral\/} if it can be represented by an $(r+1)$-tuple of integers, or,  equivalently, if $a_1-a_2$, $a_2-a_3$, \ldots, $a_r-a_{r+1}$ are all integers.   It is {\em dominant integral\/} if $a_1-a_2$, $a_2-a_3$, \ldots, $a_r-a_{r+1}$ are all non-negative integers.

\mysection{Patterns}\mylabel{s:patterns}
\noindent 
Given two non-increasing sequences $\pseq$: $\pi_1\geq\ldots\geq\pi_n$ and $\mseq$: $\mu_1\geq\ldots\geq\mu_{n-1}$ of real numbers,  of lengths $n$ and $n-1$ (for some integer $n>1$),  we say that they {\em interlace\/} if $\pi_i\geq\mu_i\geq\pi_{i+1}$ for all $i$, $1\leq i< n$.

A pattern~$\pattern$ is a sequence $\pseq^1$, \ldots, $\pseq^{r+1}$ of non-increasing sequences of real numbers, such that:
\begin{itemize}\item $\pseq^j$ has length $j$ for all $j$, $1\leq j\leq r+1$; and
		\item $\pseq^{j+1}$ and $\pseq^{j}$ interlace for all $j$, $1\leq j<r+1$.
\end{itemize}
We write $\pseq^j$ as $\pi^j_1\geq\ldots\geq\pi^j_j$, so that the interlacing condition above becomes:
\begin{equation}\label{e:pattern}
	\pi^j_i\geq \pi^{j-1}_i \quad\quad\textup{and}\quad\quad \pi^{j-1}_i\geq\pi^{j}_{i+1}\quad\quad \textup{$\forall$ $1\leq i<j\leq r+1$.}
\end{equation}

Let $\pattern$: $\pseq^1$, \ldots, $\pseq^{r+1}$ be a pattern.  Its last sequence $\pseq^{r+1}$ is called its {\em bounding sequence\/}.   We say that $\pattern$ is {\em integral\/} if 
\begin{equation}\label{e:pattern}
	\textup{the non-negative real numbers $\pi^j_i - \pi^{j-1}_i$ and $\pi^{j-1}_i-\pi^{j}_{i+1}$ are all integers for $1\leq i<j\leq {r+1}$.}
\end{equation}
Observe that the weight associated to the bounding sequence $\pseq^{r+1}$ of an integral pattern is an integral weight (in the sense of~\S\ref{ss:tuples}).    In the sequel, we will be interested exclusively in integral patterns with bounding sequence $\lseq$, where $\lseq$ is fixed as in~\S\ref{s:notation}.     We denote the set of these patterns by $\patl$.
\mysubsection{Depicting a pattern}\mylabel{ss:depictpat}
\noindent
A pattern is usually depicted as we illustrate now by means of an example.   Let $r=3$ and $\lseq=(8,6,3,1)$.    Then the sequence $4$; $5\geq4$; $7\geq5\geq2$; $8\geq6\geq3\geq1$ of sequences is a pattern belonging to~$\patl$.   It is depicted as follows:
\[
	\begin{array}{ccccccc}
		&&&4\\
		&&5&&4\\
		&7 && 5 && 2\\
		8 && 6 && 3 &&1
\end{array}
\]

Let $\pattern$: $\lseq^1$, \ldots, $\lseq^{r+1}=\lseq$ be a pattern in~$\patl$.   The sequence $\lseq^k$ is sometimes referred to as the {\em $k^\textup{th}$ row\/} of~$\pattern$.   This terminology is justified by the depiction of~$\pattern$.  It is also convenient to refer to $\lambda^k_i$ as the entry in {\em position\/} $i$ on row~$k$ of~$\pattern$ (where $\lseq^k$ is $\lambda^k_1\geq \ldots\geq \lambda^k_k$).

The {\em weight\/} of~$\pattern$ is $\mu_1\epsilon_1+\cdots+\mu_{r+1}\epsilon_{r+1}$,  where $\mu_k=(\sum_{i=1}^k \lambda^k_i) - (\sum_{i=1}^{k-1}\lambda^{k-1}_i)$.  In other words,  $\mu_k$ is the difference of the sums of the entries in the $k^\textup{th}$ row and the $(k-1)^\textup{st}$ row of~$\pattern$.   For example, the weight of the pattern depicted above is $4\epsilon_1+5\epsilon_2+5\epsilon_3+4\epsilon_4$.

As is easily seen,  the weight of a pattern~$\pattern$ in~$\patl$ is integral.    
As is also easily seen,  there is a unique pattern in~$\patl$ of weight~$\lambda$.   For example,  in case $r=4$ and $\lseq=(8,6,3,1)$,  this unique pattern is depicted below:
\[
	\begin{array}{ccccccc}
		&&&8\\
		&&8&&6\\
		&8 && 6 && 3\\
		8 && 6 && 3 &&1
\end{array}
\]

\mysection{The Gelfand-Tsetlin (GT) basis for $V(\lambda)$}\mylabel{s:gt}\noindent
Theorem~\ref{t:gt} below follows Molev's reformulation~\cite[Theorem~4.3]{molevhb} of the original theorem of Gelfand-Tsetlin~\cite{GT}.
	\begin{theorem}\mylabel{t:gt} Let $\patl$ denote the set of integral patterns with bounding sequence~$\lseq$ (see~\S\ref{s:patterns}).
	There exists a basis $\{\gt_\pattern\}$ of $V(\lambda)$ indexed by patterns~$\pattern$ in~$\patl$ such that
\begin{itemize}
\item if $\pattern$ is the unique pattern in~$\patl$ of weight $\lseq$,  then
$\gt_\pattern=v_\lambda$;
\item for every $k$, $1\leq k<r+1$, we have:
\begin{align}
	\label{e:gt1}
	(E_{k,k}-E_{k+1,k+1})\, \gtpattern&=\big(
	(\sum_{i=1}^k\lambda_i^k -\sum_{i=1}^{k-1}\lambda^{k-1}_i) -
	(\sum_{i=1}^{k+1}\lambda_i^{k+1} -\sum_{i=1}^{k}\lambda^{k}_i) 
	\big) \gtpattern, \\
	\label{e:gt2}
E_{k,k+1}\, \gtpattern
&=-\sum_{i=1}^k \frac{(\ellki-\ell^{k+1}_1)\cdots (\ellki-\ell_{k+1}^{k+1})}%
{(\ellki-\ell^k_1)\cdots\wedge\cdots(\ellki-\ell^k_k)}
	\xi_{\pattern+\delta^k_i}, \textup{and}\\
	\label{e:gt3}
E_{k+1,k}\, \gtpattern
&=\sum_{i=1}^k \frac{(\ellki-\ell^{k-1}_1)\cdots (\ellki-\ell_{k-1}^{k-1})}
{(\ellki-\ell^k_1)\cdots\wedge\cdots(\ellki-\ell^k_k)}
\gt_{\pattern-\delta_i^k}.
\end{align}
\end{itemize}
where 
\begin{itemize}
\item $E_{i,j}$ denotes the $(r+1)\times(r+1)$ matrix with its only
non-zero entry being $1$ in position $(i,j)$,
\item with $\pattern$ equal to $\lseqone$, \ldots,
$\lseqr$, $\lseqrpone$,   we define $\ell^k_i:=\lambda_i^k-i+1$,
\item $\wedge$ indicates that the vanishing factor $(\ell^k_i-\ell^k_i)$ in the denominator is
skipped,
\item $\pattern\pm\delta^k_i$ is obtained from $\pattern$ by replacing
$\lambda^k_i$ by $\lambda^k_i\pm1$,  and
\item $\gt_{\pattern\pm\delta^k_i}$ is understood to be $0$ if $\pattern\pm\delta^k_i$ is not a pattern.
\end{itemize}
\end{theorem}
	\bremark   Since $E_{k,k}-E_{k+1,k+1}$, $E_{k,k+1}$, and $E_{k+1,k}$ generate~$\lieg$ as a Lie algebra as $k$ ranges over $1\leq k< r+1$,  the equations~\eqref{e:gt1}--\eqref{e:gt3} determine completely the action of~$\lieg$ on~$V(\lambda)$.  Equation~\eqref{e:gt1} is equivalent to saying that $\gtpattern$ is a weight vector of weight equal to the weight of~$\pattern$. Thus we obtain as an immediate consequence the following formula for the character of~$\Vlambda$:
	\begin{equation}\label{e:char}
		\textup{character of~$\Vlambda$} \ = \sum_{\pattern\in\patl} \exp{(\weight(\pattern)).}
	\end{equation}
	\eremark

\mysection{The Chari-Loktev (CL) basis for $V(\lambda)$}\mylabel{s:cl}\noindent
Let $\lseq^1$, $\lseq^2$, \ldots, $\lseq^r$, $\lseq^{r+1}=\lseq$ be the rows of an integral pattern $\pattern$ with bounding sequence~$\lseq$.
As in Theorem~\ref{t:gt},  put $\ell^j_i:=\lambda^j_i-i+1$.  Let $x_{ij}^-$,  for $1\leq i\leq j\leq r$, denote the $(r+1)\times(r+1)$ complex matrix all of whose entries are zero except the one in position~$(j+1,i)$ which is~$1$.   
(Apologies for introducing yet another symbol---and a perverse one at that---for the matrix that is denoted by $E_{j+1,i}$ in Theorem~\ref{t:gt},  but it is better to preserve the notation of~\cite[Theorem~4.3]{rrv:areamax} which we will be invoking.)

For $1\leq j\leq r$,  let 
\[\rhocl{\pattern}^j:=\prod_{i=1}^j
\left(x_{ij}^-\right)^{(\ell_i^{j+1}-\ell^j_i)}\]
where, for an operator $T$ and a non-negative integer $n$,  the symbol
$T^{(n)}$ denotes the {\em divided power\/} $T^n/n!$.
The order of factors in the above product is immaterial since the matrices $x_{ij}^-$ as $i$ varies (but $j$ is fixed) 
commute with each other.  

Finally, let
\begin{equation}\label{e:clmonom}
\rhocl{\pattern}:=\rhocl{\pattern}^1\cdot\rhocl{\pattern}^2\cdot\ \cdots\ \cdot\rhocl{\pattern}^{r-1}\cdot\rhocl{\pattern}^r 
	\quad\quad\textup{and}\quad\quad
\vcl{\pattern}:=\rhocl{\pattern} \vlambda.\end{equation}
Note that the order of the factors does matter in the expression for $\rhocl{\pattern}$.
The following theorem 
is an immediate consequence of the reformulation by the present authors~\cite[Theorem~4.3]{rrv:areamax} of the original theorem of Chari-Loktev~\cite[Theorem~2.1.3]{cladv2006}.  We shall call $\rhocl{\pattern}$ the {\em Chari-Loktev (or CL)} monomial corresponding to~$\pattern$ and, in view of the theorem, $\vcl{\pattern}$ the
 {\em Chari-Loktev (or CL)} basis element corresponding to~$\pattern$.
\begin{theorem}\mylabel{t:cl}   As $\pattern$ varies over the set~$\patl$ of integral patterns with bounding sequence~$\lseq$,  the vectors $\vcl{\pattern}$ form a basis of the irreducible representation~$V(\lambda)$ of~$\lieg=\slrpone$ with highest weight $\lambda$.
\end{theorem}

\mysection{Triangular relationship between GT and CL bases of~$\Vlambda$}\mylabel{s:triangular}
\noindent
We state and prove Theorem~\ref{t:triangular} which is the goal of the present note.
	\mysubsection{Statement of the result}\mylabel{ss:tristatement} On the set~$\patl$ of integral patterns with bounding sequence~$\lseq$, we impose a partial order which we call the {\em row-wise dominance\/} order and denote by~$\geq$.   
For $\pattern$: $\lseqone$, \ldots, $\lseqr$, $\lseqrpone=\lseq$ and $\pattern'$:
$\lpseqone$, \ldots, $\lpseqr$, $\lpseqrpone=\lseq$ in $\patl$,   we let:
	\begin{quote} $\pattern'\geq\pattern$ \quad if\quad\quad $\lseq'^j\geq\lseq^j$
 \quad\quad for all\quad $1\leq j\leq r+1$, where\\
  $\lseq'^j\geq\lseq^j$ \quad if\quad\quad
$\lambda'^j_1+\cdots+\lambda'^j_i\geq
\lambda^j_1+\cdots+\lambda^j_i $ \quad\quad for all\quad
$1\leq i\leq j$.
 \end{quote}
\begin{theorem}\mylabel{t:triangular} For any pattern $\pattern$: $\lseqone$, \ldots, $\lseqr$, $\lseqrpone=\lseq$ in~$\patl$,
	in the expression for the corresponding Chari-Loktev basis element $v_\pattern$ of $V(\lambda)$ (see Theorem~\ref{t:cl}) as a linear combination of the Gelfand-Tsetlin basis $\{\gtpatternp\}_{\patternp\in\patl}$ (see Theorem~\ref{t:gt}),  only those $\gtpatternp$ with $\patternp\geq\pattern$ appear.   In other words,  we can write: 
\begin{equation}\label{e:triangular}
\vcl{\pattern}=\sum_{\pattern'\geq\pattern} \coeff^{\pattern'}_{\pattern}\gt_{\pattern'}
\end{equation}
	where $\coeff^{\patternp}_\pattern$ are entries of the transition matrix between the two bases.  The diagonal element $\coeff^\pattern_\pattern$ in the transition matrix is given by:
\begin{equation}\label{e:tricoeff}
\coeff^\pattern_\pattern = \prod_{1\leq i<j\leq
r+1}
\prod_{p=i}^{j-1}\frac{(\ell^{j-1}_i-\ell^{j-1}_p)!}%
{(\ell^j_i-\ell_{p}^{j-1})!}
\end{equation}
	where,  as in Theorem~\ref{t:gt},  $\ell^j_i:=\lambda^j_i-i+1$.
\end{theorem}
\mysubsection{Proof of Theorem~\ref{t:triangular}}\mylabel{ss:pf:triangular}
Proceed by induction on the length of $\pattern$,  where length is
defined by 
\begin{equation}\label{e:length}
\lengthp{\pattern}:= \sum_{1\leq i<j\leq r+1}(\lambda^j_i-\lambda^{j-1}_i)
= \sum_{1\leq i<j\leq r+1}(\ell^j_i-\ell^{j-1}_i).
\end{equation}
If the length is $0$,   then $\pattern$ is the unique one of 
weight $\lseq$,  the only pattern $\pattern'$ such
that $\pattern'\geq\pattern$ is $\pattern$ itself,
$\vcl{\pattern}=\gtpattern=v_\lambda$, and $\coeff^\pattern_\pattern=1$, and we
are done.

\begin{figure}[h]
\begin{center}
 \setlength{\unitlength}{0.0007489in}
\begingroup\makeatletter\ifx\SetFigFont\undefined%
\gdef\SetFigFont#1#2#3#4#5{%
  \reset@font\fontsize{#1}{#2pt}%
  \fontfamily{#3}\fontseries{#4}\fontshape{#5}%
  \selectfont}%
\fi\endgroup%
{\renewcommand{\dashlinestretch}{30}
\begin{picture}(7587,4726)(-1500,10)
\drawline(7575,1864)(7575,1864)
\drawline(7575,1864)(7575,1864)
\drawline(3210,4699)(3210,4699)
\drawline(3210,4699)(3210,4699)
\drawline(60,379)(60,379)
\drawline(960,379)(2490,3079)
\drawline(1005,379)(2535,3034)
\drawline(2895,2179)(2895,2179)
\drawline(3480,1279)(3480,1279)
\drawline(2940,2179)(2940,2179)
\drawline(2895,2179)(2895,2179)
\drawline(2985,2179)(2985,2179)
\drawline(2265,829)(3030,2179)
\drawline(3390,1234)(3390,1234)
\drawline(3210,1729)(3210,1729)
\drawline(3255,1729)(3255,1729)
\drawline(3480,1324)(3480,1324)
\drawline(3255,1729)(2715,829)
\drawline(3480,1279)(3255,829)
\drawline(3750,829)(3750,829)
\drawline(3435,1279)(3435,1279)
\drawline(3210,1729)(3210,1729)
\drawline(3300,829)(3525,1279)
\drawline(2760,829)(3210,1594)
\drawline(3210,1594)(3255,1684)
\drawline(2715,2629)(2715,2629)
\drawline(2985,2179)(2985,2179)
\drawline(1185,19)(2850,2854)
\drawline(2220,829)(2985,2179)
\drawline(60,379)(60,379)
\drawline(3750,829)(510,829)
\drawline(1410,2629)(2760,2629)
\drawline(1545,2179)(2985,2179)
\drawline(1320,1729)(3300,1729)
\drawline(1005,1279)(3525,1279)
\drawline(2760,2629)(2760,2629)
\drawline(1725,829)(2760,2629)
\drawline(5010,3574)(5010,3574)
\drawline(2085,3079)(2535,3079)
\drawline(4065,379)(2265,3484)(2265,3439)
	(2265,3484)(4065,379)(2130,379)
	(825,379)(195,379)
\drawline(510,379)(2265,3394)
\drawline(2265,3484)(465,379)
\put(2895,2854){\makebox(0,0)[lb]{\smash{{\SetFigFont{12}{14.4}{\rmdefault}{\mddefault}{\updefault}$i_0$}}}}
\put(1275,2584){\makebox(0,0)[lb]{\smash{{\SetFigFont{12}{14.4}{\rmdefault}{\mddefault}{\updefault}$i_0$}}}}
\put(15,739){\makebox(0,0)[lb]{\smash{{\SetFigFont{12}{14.4}{\rmdefault}{\mddefault}{\updefault}$j_0-1$}}}}
\put(15,334){\makebox(0,0)[lb]{\smash{{\SetFigFont{12}{14.4}{\rmdefault}{\mddefault}{\updefault}$j_0$}}}}
\end{picture}
}
\end{center}
\caption{The figure shows the first $j_0$ rows of the
pattern~$\pattern$.
The Southwest--Northeast double lines indicate equalities, while the lone single Southwest--Northeast line between rows $j_0$ and $j_0-1$ indicates $\gneq$.}\label{f:pattern}
\end{figure}
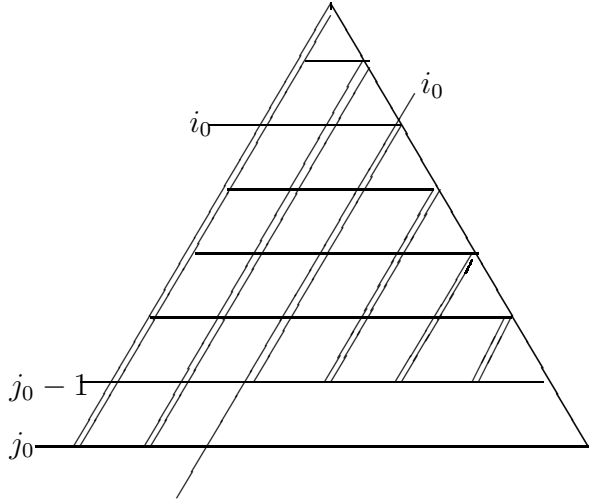

Let now the length of $\pattern$ be at least one.    Let
$j_0:=\min\{j\geq 2\st \lambda^j_i\neq\lambda^{j-1}_i\textup{ for some
$1\leq i\leq j$}\}$ and
$i_0:=\min\{i\st \lambda_i^{j_0}\neq\lambda_i^{j_0-1}\}$.  Let
$\tpattern$: $\tlseq^1$, \ldots, $\tlseq^r$, $\tlseq^{r+1}=\tlseq$ be
obtained from $\pattern$ by replacing $\lambda^j_{i_0}$ by
$\lambda^j_{i_0}+1$ for $i_0\leq j< j_0$ (the other entries stay
the same).    See Figure~\ref{f:pattern}.

It is easily seen that
\[
\textup{ $\tpattern$ is a pattern, \quad
$\tpattern\geq\pattern$,\quad and\quad $\lengthp{\tpattern}=\lengthp{\pattern}-1$.}
\]
It is also easily seen from the definitions in~\S\ref{s:cl} that the
Chari-Loktev monomials of $\pattern$ and $\tpattern$, and consequently
also the corresponding basis elements, are related thus:
\begin{equation}\label{e:trpf:1}
	\rhocl{\pattern}=\frac{1}{\ell^{j_0}_{i_0}-\ell^{j_0-1}_{i_0}}x_{i_0,j_0-1}^-\rhocl{\tpattern}\quad\quad\textup{and}\quad\quad
\vcl{\pattern}=\frac{1}{\ell^{j_0}_{i_0}-\ell^{j_0-1}_{i_0}}E_{j_0,i_0}\vcl{\tpattern}.
\end{equation}

By the induction hypothesis, we have:
\begin{gather}\label{e:trpf:2}
\vcl{\tpattern}=\sum_{\pattern'\geq\tpattern}\coeff_{\tpattern}^{\pattern'}\gt_{\pattern'}
	\quad\quad\textup{and}\quad\quad
\coeff^{\tpattern}_{\tpattern} = 
\prod_{1\leq i<j\leq r+1}
\prod_{p=i}^{j-1}\frac{(\tell^{j-1}_i-\tell^{j-1}_p)!}%
{(\tell^j_i-\tell_{p}^{j-1})!}
\end{gather}
where $\tell^j_i:=\tlambda^j_i-i+1$.   It is convenient to rewrite the expression above for
$\coeff^{\tpattern}_{\tpattern}$ in terms of the numbers $\ell^j_i=\lambda^j_i-i+1$
for $\pattern$.       We claim that
\begin{equation}\label{e:trpf:coeff:00}
\coeff^{\tpattern}_{\tpattern}=
(\ell^{j_0}_{i_0}-\ell^{j_0-1}_{i_0})
(\prod_{p=i_0+1}^{j_0-1}(\ell_{i_0}^{j_0-1}-\ell^{j_0-1}_p+1))
\left(\prod_{1\leq i<j\leq r+1}\prod_{p=i}^{j-1}\frac{(\ell_i^{j-1}-\ell^{j-1}_p)!}%
{(\ell^j_i-\ell_{p}^{j-1})!}\right).
\end{equation}
To prove the claim, we first observe the following:
\begin{gather}\label{e:trpf:coeff:0}
\prod_{p=i}^{j-1}\frac{(\tell_i^{j-1}-\tell^{j-1}_p)!}%
{(\tell^j_i-\tell_{p}^{j-1})!}
\quad=\quad
\prod_{p=i}^{j-1}\frac{(\ell_i^{j-1}-\ell^{j-1}_p)!}%
{(\ell^j_i-\ell_{p}^{j-1})!} \quad\quad\textup{for
$(j,i)\neq(j_0,i_0)$.}
\end{gather}
Indeed, for $j\geq j_0$,   we have $\lseq^j=\tlseq^j$,  so (\ref{e:trpf:coeff:0})
holds for $j>j_0$;    it also holds for $i>i_0$
since $\ell^j_i=\tell^j_i$ in that case;  finally,  if
either $j<j_0$ or $j=j_0$ and $i<i_0$,   then
$\ell^j_i=\ell^{j-1}_i$
and $\tell^j_i=\tell^{j-1}_i$,  so the product on either side of (\ref{e:trpf:coeff:0})
 equals $1$. 

	Continuing with the proof  of  the claim
 (\ref{e:trpf:coeff:00}),    we calculate the left hand side of
 (\ref{e:trpf:coeff:0}) in the case $(j,i)=(j_0,i_0)$:
\begin{align*}
&\prod_{p=i_0}^{j_0-1}\frac{(\tell_{i_0}^{j_0-1}-\tell^{j_0-1}_p)!}%
{(\tell^{j_0}_{i_0}-\tell_{p}^{j_0-1})!} =
\frac{1}{(\tell^{j_0}_{i_0}-\tell^{j_0-1}_{i_0})!}\times
\prod_{p=i_0+1}^{j_0-1}\frac{(\tell_{i_0}^{j_0-1}-\tell^{j_0-1}_p)!}%
{(\tell^{j_0}_{i_0}-\tell_{p}^{j_0-1})!}\\
 &=
\frac{1}{(\ell^{j_0}_{i_0}-\ell^{j_0-1}_{i_0}-1)!}\times
\prod_{p=i_0+1}^{j_0-1}\frac{(\ell_{i_0}^{j_0-1}-\ell^{j_0-1}_p+1)!}%
{(\ell^{j_0}_{i_0}-\ell_{p}^{j_0-1})!}\\
&=
\frac{(\ell^{j_0}_{i_0}-\ell^{j_0-1}_{i_0})}%
{(\ell^{j_0}_{i_0}-\ell^{j_0-1}_{i_0})!}\times
\left(\prod_{p=i_0+1}^{j_0-1}(\ell_{i_0}^{j_0-1}-\ell^{j_0-1}_p+1)\right)
\times
\left(
\prod_{p=i_0+1}^{j_0-1}\frac{(\ell_{i_0}^{j_0-1}-\ell^{j_0-1}_p)!}%
{(\ell^{j_0}_{i_0}-\ell_{p}^{j_0-1})!}  \right)\\
&=
(\ell^{j_0}_{i_0}-\ell^{j_0-1}_{i_0})
\times
\left(\prod_{p=i_0+1}^{j_0-1}(\ell_{i_0}^{j_0-1}-\ell^{j_0-1}_p+1)\right)
\times
\left(
\prod_{p=i_0}^{j_0-1}\frac{(\ell_{i_0}^{j_0-1}-\ell^{j_0-1}_p)!}%
{(\ell^{j_0}_{i_0}-\ell_{p}^{j_0-1})!}  \right)\\
\end{align*}
where, to justify the last equality, we observe that
 the $(\ell^{j_0}_{i_0}-\ell^{j_0-1}_{i_0})!$ in the
denominator of the first factor is incorporated into the last factor
by changing the lower range of $p$ from $i_0+1$ to $i_0$.   This
calculation along with (\ref{e:trpf:coeff:0}) proves the
claim~(\ref{e:trpf:coeff:00}).

We now  prove (\ref{e:triangular}) of the theorem.
Substituting the expression for $\vcl{\tpattern}$ in~(\ref{e:trpf:2})
into the expression for $\vcl{\pattern}$ in~(\ref{e:trpf:1}),  we get
\begin{equation}\label{e:trpf:3}
\vcl{\pattern} =\frac{1}{\ell^{j_0}_{i_0}-\ell_{i_0}^{j_0-1}}
\sum_{\pattern'\geq\tpattern} \coeff^{\pattern'}_{\tpattern}E_{j_0,i_0}\gt_{\pattern'}.
\end{equation}
Observe that
\begin{equation}\label{e:lieprod}E_{j_0,i_0}=[E_{j_0,j_0-1}[E_{j_0-1,j_0-2},\cdots [E_{i_0+2,i_0+1},E_{i_0+1,i_0}]\cdots]].
\end{equation}
So the operator~$E_{j_0,i_0}$ on~$V(\lambda)$ is a linear combination with coefficients $\pm1$ of the operators
\begin{equation}\label{e:eprod}
	E_{\sigma(j_0),\sigma(j_0)-1}E_{\sigma(j_0-1),\sigma(j_0-1)-1}\cdots E_{\sigma(i_0+1),\sigma(i_0+1)-1}
\end{equation}
as $\sigma$ varies over permutations of the set $\{i_0+1,\ldots,j_0\}$,  with $1$ being the coefficient of the operator \[E_{j_0,j_0-1}E_{j_0-1,j_0-2}\cdots E_{i_0+2,i_0+1}E_{i_0+1,i_0}.\]
Now, \eqref{e:triangular} follows from formula~\eqref{e:gt3} for the action of $E_{k+1,k}$ in Theorem~\ref{t:gt} and item~\eqref{i:trione} of the following lemma.
\begin{lemma}\mylabel{l:triangular}
Let $\pattern'$ be any pattern in~$\patl$ 
such that $\pattern'\geq\tpattern$.
Let $\pattern''$ be any pattern obtained from $\pattern'$ by
decreasing by $1$ exactly one entry in each of the rows $i_0$, \ldots,
$j_0-1$.   Then
	\begin{enumerate}\item\label{i:tri0}$\pattern''\in\patl$. \item\label{i:trione}
			$\pattern''\geq\pattern$.   
		\item\label{i:tritwo} 
$\pattern''=\pattern$ only if $\pattern'=\tpattern$.
\item\label{i:trithree} The coefficient of $\gtpattern$ in~$E_{j_0,i_0}\gtpatternt$ equals that of $\gtpattern$ in 
	$E_{j_0,j_0-1}E_{j_0-1,j_0-2}\cdots E_{i_0+1,i_0}\gtpatternt$,  and consequently:
\begin{equation}\label{e:coeffexp}
	\textup{coefficient of~$\gtpattern$ in~$E_{j_0,i_0}\gtpatternt$}\ =\ 
	\prod_{k=i_0}^{j_0-1}\textup{coefficient of $\gt_{\tpattern(k+1)}$ in $E_{k+1,k}\,\gt_{\patternt(k)}$}
\end{equation}
			where, for $k$ in the range $i_0\leq k\leq j_0$, $\tpattern(k)$ denotes the pattern whose initial rows up to but not including the $k^\textup{th}$ are those of~$\pattern$ and the remaining ones are those of~$\tpattern$ (so that $\tpattern(i_0)=\tpattern$ and $\tpattern(j_0)=\pattern$).
	\end{enumerate}
\end{lemma}
\begin{myproof}
	Since $j_0-1<r+1$,  it is clear that $\pattern''$ shares its $(r+1)^\textup{st}$ row with~$\pattern'$ and so~\eqref{i:tri0} is clear.
Let the rows of $\pattern'$ be $\lseq'^1$, \ldots, $\lseq'^r$,
$\lseq'^{r+1}=\lseq$
and those of $\pattern''$ be $\lseq''^1$, \ldots, $\lseq''^r$,
$\lseq''^{r+1}=\lseq$.

	Towards the proof of~\eqref{i:trione}, first observe that for $j$ such that $1\leq j<i_0$ or $j_0\leq j\leq r+1$,  we have 
\[\lseq''^j=\lseq'^j\geq\ltildeseq^j=\lseq^j.
\]
It remains to show that for $j$ such that $i_0\leq j<j_0$ and $k$ such
that $1\leq k\leq j$,   it holds that
\begin{equation}\label{e:trpf''3}
\lambda''^j_1+\cdots+\lambda''^j_k\geq \lambda^j_1+\cdots+\lambda^j_k.
\end{equation}
For such $j$,  let $r_j$ denote the position in row $j$ in which $\pattern''$ differs from $\pattern'$.    

We consider three cases:   $k<r_j$,  $i_0\leq k$,  and $r_j\leq k<i_0$
(this last may not occur if $i_0\leq r_j$).    In the first case we
have 
\[
 \lambda''^j_1+\cdots+\lambda''^j_k
= \lambda'^j_1+\cdots+\lambda'^j_k
\geq \tlambda^j_1+\cdots+\tlambda^j_k
\geq \lambda^j_1+\cdots+\lambda^j_k
\]
and in the second case
\[
 \lambda''^j_1+\cdots+\lambda''^j_k
\geq \lambda'^j_1+\cdots+\lambda'^j_k-1
\geq \tlambda^j_1+\cdots+\tlambda^j_k-1
=\lambda^j_1+\cdots+\lambda^j_k.
\]
In the third case,  using the fact that $\pattern''$ is a pattern, we have:
\[
 \lambda''^j_1+\cdots+\lambda''^j_k
\geq \lambda''^k_1+\cdots+\lambda''^k_k
= \lambda'^k_1+\cdots+\lambda'^k_k
\geq \tlambda^k_1+\cdots+\tlambda^k_k
=\lambda^k_1+\cdots+\lambda^k_k
=\lambda^j_1+\cdots+\lambda^j_k
\]
	which completes the proof of~\eqref{i:trione}. 

	Let us now prove~\eqref{i:tritwo}.    Since $\pattern''$ is given to be $\pattern$,   it follows that $\pattern'$ differs from $\pattern$ in exactly one position on each of the rows $i_0$, \ldots, $j_0-1$ and furthermore the difference in each case is only that the entry of $\pattern'$ is precisely one more than the corresponding entry of~$\pattern$.    Thus,  as far as other rows are concerned, that is, for $j$ not in the range $i_0\leq j<j$,   we have $\lseq'^j=\lseq^j=\tlseq^j$.    

	Now let $j$ be fixed such that $i_0\leq j<j$.   We need to show $\lseq'^j=\tlseq^j$, or, in more detail, $\lambda'^j_i = {\tlambda}_i^j$ for all $i$, $1\leq i\leq j$.  Suppose that we have proved this for $i$ in the range $1\leq i< i_0$.     Then, since
	\[\lambda'^j_1+\cdots+\lambda'^j_{i_0-1}+\lambda'^j_{i_0}\geq \tlambda^j_1+\cdots+\tlambda^j_{i_0-1}+\lambda^j_{i_0}\]
	(this equation holds because $\pattern'\geq\tpattern$),   it would follow that $\lambda'^j_{i_0}\geq\tlambda^j_{i_0}=\lambda^j_{i_0}+1$.   But given how little $\pattern'$ is allowed to differ from $\pattern$,   it would follow that equality holds here,  that is, $\lambda'^j_{i_0}=\tlambda^j_{i_0}=\lambda^j_{i_0}+1$.   Moreover it would follow for the same reason that $\lambda'^j_i=\lambda^j_i=\tlambda^j_i$ for $i$ in the range $i>i_0$ as well, and we would be done with the proof of~\eqref{i:tritwo}.

	Thus it only remains to prove that $\lambda'^j_i=\tlambda^j_i$ for $i$ in the range $1\leq i<i_0$,   which we proceed to do by induction on~$i$.   The induction hypothesis is that $\lambda'^j_s=\tlambda^j_s$ for $1\leq s<i$ (which is vacuously true in the base case $i=1$).    Since 
	\[\lambda'^j_1+\cdots+\lambda'^j_{i-1}+\lambda'^j_{i}\geq \tlambda^j_1+\cdots+\tlambda^j_{i-1}+\lambda^j_{i}\]
	(this equation holds because $\pattern'\geq\tpattern$),   it follows that $\lambda'^j_{i}\geq\tlambda^j_{i}$.   Using this (as the fourth (in)equality in the string of (in)equalities below), we have:
\[\lambda^{j_0}_i=\lambda''^{j_0}_i=\lambda'^{j_0}_i\geq\lambda'^j_i \geq\tlambda^j_i=\lambda^j_i=\lambda^{j_0}_i\quad\quad\textup{so that $\lambda'^j_i=\tlambda_i^j(=\lambda^j_i=\lambda^{j_0}_i)$}\]
	where, the first equality is justified because $\pattern''$ is given to be equal to~$\pattern$;  the second because $\pattern''$ and $\pattern'$ have the same $j_0^\textup{th}$ row (among others);  the third because $\pattern'$ is a pattern;  the fifth because $\tpattern$ differs from $\pattern$ only in positions $i_0$ in rows $i_0$, \ldots, $j_0-1$;   and, finally, the sixth by the choice of $i_0$ and $j_0$.     This finishes the proof of~\eqref{i:tritwo}.

	Turning to~\eqref{i:trithree},   given how $E_{j_0,i_0}$ is a linear combination of operators as in~\eqref{e:eprod} and the formula~\eqref{e:gt3} for the action of~$E_{k+1,k}$ on a Gelfand-Tsetlin basis vector,  it suffices to prove the following:  
	\begin{quote}
		Suppose that $\patternq_{i_0}$, \ldots, $\patternq_{j_0}$ be a sequence of patterns (indexed by $k$ in the range $i_0\leq k\leq j_0$) and $\sigma$ a permutation of the set~$\{i_0,\ldots,j_0-1\}$ such that $\patternq_{i_0}=\tpattern$, $\patternq_{j_0}=\pattern$,  and, for each $k$, $i_0\leq k<j_0$, the pattern~$\patternq_{k+1}$ is obtained from $\patternq_k$ by making only one change, namely, decreasing the entry in position~$i_0$ in the row~$\sigma(k)$.   Then $\sigma$ is the identity permutation.
	\end{quote}
	By way of contradiction,  suppose that $\sigma$ is not the identity.	Let $k$ be least such that $\sigma(k)\neq k$.   Then $\sigma(k)>k$.      Compare the entries in position $i_0$ on rows $\sigma(k)$ and $k$ in the pattern~$\patternq_{k+1}$.     The former is $1$ less than the latter, which contradicts the assumption that $\patternq_{k+1}$ is a pattern.   This finishes the proof of~\eqref{i:trithree} and so also of the lemma.
	\emyproof

We now turn to the proof of (\ref{e:tricoeff}) of
Theorem~\ref{t:triangular}.   Combining item~\eqref{i:tritwo} of Lemma~\ref{l:triangular} with
(\ref{e:trpf:3}),   we obtain the following:
\begin{equation}\label{e:trpf:coeff:4}
\coeff^\pattern_\pattern=\frac{\coeff^{\tpattern}_{\tpattern}}{\ell^{j_0}_{i_0}-\ell^{j_0-1}_{i_0}}\times\textup{the
coefficient of $\gt_\pattern$ in $E_{j_0,i_0}\gt_{\tpattern}$.}
\end{equation}
Substituting for $\coeff^{\tpattern}_{\tpattern}$ its value from (\ref{e:trpf:coeff:00}),  we see that it suffices
to prove the following:
\begin{lemma}\mylabel{l:trpf:2}  The coefficient of $\gt_{\pattern}$
in $E_{j_0,i_0}\gt_{\tpattern}$ equals
\begin{equation}
\prod_{p=i_0+1}^{j_0-1}
\frac{1}{
(\ell_{i_0}^{j_0-1}-\ell^{j_0-1}_p+1)}.
\end{equation}
\end{lemma}
\bmyproof
Item~\eqref{i:trithree} of Lemma~\ref{l:triangular} says that the required coefficient equals
\[
	\prod_{p=i_0}^{j_0-1} \textup{coefficient of $\gt_{\tpattern(p+1)}$ in $E_{p+1,p}\,\gt_{\tpattern(p)}$}.
	\]
From~\eqref{e:gt3} in Theorem~\ref{t:gt},  we see that, for $p$ in the range $i_0\leq p<j_0$,
the coefficient of $\gt_{\tpattern(p+1)}$ in $E_{p+1,p}\,\gt_{\tpattern(p)}$ is
\[ 
	\frac{\prod_{q=1}^{p-1}(\alpha_{i_0}+1-\alpha_{q})}%
	{\prod_{q=1}^{i_0-1}(\alpha_{i_0}+1-\alpha_q) \times \prod_{q=i_0+1}^p(\alpha_{i_0}+1-\alpha_q)}
	=\frac{1}{(\alpha_{i_0}+1-\alpha_p)} =\frac{1}{(\ell^{j_0-1}_{i_0}-\ell^{j_0-1}_{p}+1)}
\] 
where we have denoted by $\alpha_q$ the value $\ell^q_q=\ell^{q+1}_q=\ldots=\ell^{j_0-1}_q$ (for $q$ in the range $1\leq q<j_0$).   This finishes the proof of Lemma~\ref{l:trpf:2} and so also that of Theorem~\ref{t:triangular}.
\emyproof

\mysection{An example illustrating Theorem~\ref{t:triangular}}\mylabel{s:example}
\noindent  Let $\lieg=sl_3$ (that is, $r=2$).  Let $\lambda=4\epsilon_1+2\epsilon_2$ and $\lseq=4\geq2\geq0$.   In $\patl$ there are three patterns of weight $0$ and here they are listed in row-dominance order:
\[
\begin{array}{ccccc}
& & 2 \\
& 4 & & 0\\
4 &&2 &&0
\end{array}
\quad\geq\quad
\begin{array}{ccccc}
& & 2 \\
& 3 & & 1\\
4 &&2 &&0
\end{array}
\quad\geq\quad
\begin{array}{ccccc}
& & 2 \\
& 2 & & 2\\
4 &&2 &&0
\end{array}
\]
Let us denote them by $\pattern_1$, $\pattern_2$, and $\pattern_3$ respectively.   Letting $\gt_1$, $\gt_2$, and $\gt_3$ be the corresponding GT basis elements of~$\Vlambda$ and $\vcl{1}$, $\vcl{2}$, $\vcl{3}$ the corresponding Chari-Loktev basis elements,   we have:
\[
	\left(
	\begin{array}{c}
		\vcl{1}\\
		\vcl{2}\\
		\vcl{3}
	\end{array}
	\right)
	=
	\left(
	\begin{array}{rrr}
	1/4 & 0 & 0 \\
	-1/4 & 1/4 & 0 \\
	1/24 & -1/8 & 1/24 \\
	\end{array}
	\right)
	\left(
	\begin{array}{c}
		\gt_1\\
		\gt_2\\
		\gt_3
	\end{array}
	\right)
	\]


\bibliographystyle{bibsty-final-no-issn-isbn}
\addcontentsline{toc}{section}{References}
\ifthenelse{\equal{\finalized}{no}}{
	\bibliography{abbrev,references}}{

\begin{thebibliography}{99}
\footnotesize\itemsep=0pt
\providecommand{\eprint}[2][]{\href{http://arxiv.org/abs/#2}{arXiv:#2}}

\bibitem{cladv2006}
V.~Chari and S.~Loktev, {\em Weyl, {D}emazure and fusion modules for the
  current algebra of {$\germ s\germ l_{r+1}$}\/}, Adv. Math., {\bf 207}, no.~2,
  2006, pp.~928--960,
  \url{http://dx.doi.org/10.1016/j.aim.2006.01.012},
\href{http://arxiv.org/abs/math/0502165}{arXiv:math/0502165}.



\bibitem{GT}
I.~M. Gel{\cprime}fand and M.~L. Cetlin, {\em Finite-dimensional
  representations of the group of unimodular matrices\/}, Doklady Akad. Nauk
  SSSR (N.S.), {\bf 71},  1950, pp.~825--828.
		English transl. in: I. M. Gelfand, ``Collected papers". Vol II, Berlin: Springer-Verlag 1988, pp. 653--656.

  
\bibitem{molevhb}
A.~I. Molev, {\em Gelfand-{T}setlin bases for classical {L}ie algebras\/}, in:
  {\em Handbook of algebra. {V}ol. 4\/}, vol.~4 of Handb. Algebr.,
  Elsevier/North-Holland, Amsterdam, 2006, pp. 109--170,
  \url{http://dx.doi.org/10.1016/S1570-7954(06)80006-9},
		\href{http://arxiv.org/abs/math/0211289}{arXiv:math/0211289}. 

\bibitem{rrv:areamax}
K.~N. Raghavan, B.~Ravinder, and S.~Viswanath, {\em On {C}hari-{L}oktev bases
  for local {W}eyl modules in type {$A$}\/}, J. Combin. Theory Ser. A, {\bf
  154}, 2018, pp.~77--113,
  \url{https://doi.org/10.1016/j.jcta.2017.08.011},
		\href{http://arxiv.org/abs/1606.01191}{arXiv:math/1606.01191}.
\end{thebibliography}

}

\end{document}